\documentclass[12pt]{amsart}

\theoremstyle{plain}
  \newtheorem{theorem}{Theorem}[section]
  \newtheorem{proposition}[theorem]{Proposition}
  \newtheorem{lemma}[theorem]{Lemma}
  \newtheorem{corollary}[theorem]{Corollary}
  
\theoremstyle{definition}

\numberwithin{equation}{section}

\def\Pfaff{{\rm P}}

\def\pol{{\mathrm{pol}}}
\def\iin{{\mathrm{in}}}
\def\Hilb{{\mathrm{Hilb}}}
\def\xx{{\underline{\bf x}}}
\def\Amatrix{{A}}
\newcommand{\varx}[2]{{X}_{#1#2}}

\def\Determ{{\rm D}}

\begin{document}

\title[A spherical initial ideal for Pfaffians]{A spherical initial ideal for Pfaffians}

\author{Jakob Jonsson}

\address{
Technische Universit\"at Berlin, 
Institut f\"ur Mathematik, MA 6-2, 
Diskrete Geometrie, 
Stra\ss e des 17. Juni 136, 
10623 Berlin, Germany}
\address{ 
Department of Mathematics, 
Royal Institute of Technology,
S-10044 Stockholm, Sweden} 

\email{jakobj@kth.se}

\author{Volkmar Welker}
\address{Fachbereich Mathematik und Informatik\\
Philipps-Universit\"at Marburg\\
35032 Marburg, Germany}
\email{welker@mathematik.uni-marburg.de}

\thanks{First author supported by Graduiertenkolleg `Combinatorics, Geometry, Computation`, DFG-GRK 588/2.  
Both authors were supported by EU
      Research Training Network 
``Algebraic Combinatorics in Europe'', grant HPRN-CT-2001-00272.}

\keywords{Pfaffian ideal, Gr\"obner basis, simplicial complex, Associahedron}

\subjclass{ }

\begin{abstract}
  We determine a term order on the monomials in the variables 
  $\varx{i}{j}$, $1 \leq i < j \leq n$, such that corresponding initial 
  ideal of the ideal of Pfaffians of degree $r$ of a generic $n$ by $n$ 
  skew-symmetric matrix is the Stanley-Reisner ideal of 
  a join of a simplicial sphere and a simplex.
  Moreover, we demonstrate that the Pfaffians of the $2r$ by $2r$
  skew-symmetric submatrices form a Gr\"obner basis for the given
  term order.
  The same methods and similar term orders as for the Pfaffians also 
  yield squarefree initial ideals for certain determinantal ideals. 
  Yet, in contrast to the case of Pfaffians, the corresponding
  simplicial complexes are balls that do not decompose into a join as
  above.
\end{abstract}

\maketitle

%
%

\section{Introduction}
\label{introduction}

The ideal $\Pfaff_{n,r}$ of Pfaffians of degree $r$ of a generic $n \times n$
skew-symmetric matrix is one of the classical ideals in commutative 
algebra. Many of its properties are well understood. In particular,
Herzog and Trung \cite{HT} have constructed a term order for which the
standard generators of $\Pfaff_{n,r}$ constitute a Gr\"obner
basis. Indeed, this term order has many nice properties, including the
property that the corresponding initial ideal is squarefree and is
the Stanley-Reisner ideal of a simplicial ball. In this paper, we
determine a different term order $\preceq$ leading to the same
Gr\"obner basis and a corresponding initial ideal that is squarefree
and the Stanley-Reisner ideal of a simplicial ball. Yet, in  
contrast to Herzog and Trung's situation, our term order $\preceq$
satisfies the following conditions:

\begin{itemize} 
  \item[(S)] The initial ideal $\iin_\preceq (\Pfaff_{n,r})$ is the 
             Stanley-Reisner ideal of a simplicial complex that
             decomposes into a join $\Delta * 2^\Theta$ of a
             simplicial sphere $\Delta$ and a full simplex $2^\Theta$
             on ground set $\Theta$.
  \item[(M)] The cardinality of the set $\Theta$ is the absolute value
             of the $a$-invariant of the quotient of the polynomial 
             ring by $\Pfaff_{n,r}$.
\end{itemize}

If condition (S) is satisfied, then we say that $\preceq$ is a term
order with {\em spherical} initial ideal. 
Note that if the set $\Theta$ from condition (S) is nonempty, then 
the simplicial complex $\Delta * 2^\Theta$ is a ball.  
As mentioned above, the initial ideal of Herzog and Trung's term order
\cite{HT} is also the Stanley-Reisner ideal of a simplicial 
ball. However, this ball does not decompose into a join as postulated
in (S). 

Condition (M) assures a certain minimality of
the initial ideal. Indeed, among all squarefree ideals occurring as 
initial ideals of $\Pfaff_{n,r}$, an ideal satisfying this condition
is one that uses the least number of variables in its minimal
generating set.

Recently, term-orders satisfying (S) and (M) have been constructed 
for several classes of ideals 
\cite{A,CHT,OH,RW}; see in particular \cite{CHT} for a survey.
In all these cases, the ideals were the
defining ideals of affine semigroup rings or were reduced to this case
by a flat deformation.
This paper provides the first instance of a 
non affine semigroup ring for which this construction can be performed.

In addition to the Pfaffian case, in Section \ref{determinantal} 
we study similar term orders for ideals generated by $r \times r$ minors 
of matrices whose generic entries form a stack polyomino shape. 
Again the initial ideals are squarefree, 
but this time the corresponding simplicial complexes do not decompose
into a join of a simplicial sphere and a simplex. 

\section{Basic Facts and the Main Result}
\label{definitions}

Let $S_n = k[\varx{i}{j} : 1 \leq i < j \leq n]$ be the polynomial 
ring over the field $k$ in 
the entries of 
a generic $n \times n$ skew-symmetric matrix 
$$A_n = 
\left( 
\begin{array}{cccccc} 
0         & \varx{1}{2}   & \varx{1}{3}    & \cdots & \varx{1}{n-1}  & \varx{1}{n} \\
-\varx{1}{2}   & 0        & \varx{2}{3}    & \cdots & \varx{2}{n-1}  & \varx{2}{n} \\               
-\varx{1}{3}   & -\varx{2}{3}  &      0    & \cdots & \varx{3}{n-1}  & \varx{3}{n} \\
\cdots    & \cdots   & \cdots    & \cdots & \cdot     & \cdot \\
\cdots    & \cdots   & \cdots    & \cdots & \cdot     & \cdot \\
\cdots    & \cdots   & \cdots    & \cdots & \cdot     & \cdot \\
-\varx{1}{n-1} & -\varx{2}{n-1}& -\varx{3}{n-1} & \cdots &         0 & \varx{n-1}{n} \\
-\varx{1}{n}   &  -\varx{2}{n} & -\varx{3}{n}   & \cdots & -\varx{n-1}{n} & 0  
\end{array}
\right).$$

It is a well-known fact from linear algebra that $\det(A_n) = 0$ if
$n$ is odd and $\det(A_n) = p_{A_n}^2$ if $n$ is even for a certain
polynomial $p_{A_n} \in S_n$ of degree $n$.
The polynomial $p_{A_n}$ is called the {\em Pfaffian} of $A_n$. 
For $1 \leq r \leq \lfloor \frac{n}{2} \rfloor$ and indices 
$1 \leq j_1 < \cdots < j_{2r} \leq n$, we denote by 
$A_{j_1, \ldots, j_{2r}}$ the submatrix of $A_n$ constructed by
selecting the rows and columns indexed by $j_1, \ldots, j_{2r}$. 
Clearly, the matrices $A_{j_1, \ldots, j_{2r}}$ are the only
skew-symmetric submatrices of $A_n$ of size $2r$. We denote by
$\Pfaff_{n,r}$ the ideal in $S_n$ generated by the
Pfaffians of the matrices $A_{j_1, \ldots, j_{2r}}$, 
$1 \leq j_1 < \cdots < j_{2r} \leq n$. 

Before proceeding, we briefly mention some well-known definitions and
facts about simplicial complexes and Gr\"obner bases. We refer to 
Bruns and Herzog \cite[Section~5]{BH} for more details
on simplicial complexes and to Fr\"oberg \cite{F2} for more details on
Gr\"obner bases. 

Let $T = k[x_1, \ldots, x_\ell]$ be the polynomial ring
in $\ell$ variables $x_1, \ldots, x_\ell$. A {\em monomial} in $T$ is
a product $m = \displaystyle{\prod_{1 =1}^\ell x_{i}^{a_i}}$, where 
each $a_i$ is a nonnegative integer. The monomial $m$ is called
{\em squarefree} if $a_{i} \in \{0,1\}$ for $1 \leq i \leq \ell$. An
ideal $I$ in $T$ is called a {\em (squarefree) monomial ideal} if $I$
is generated by a set of (squarefree) monomials. 
The {\em Stanley-Reisner ideal} $I_\Delta \subseteq T$ of a
simplicial complex $\Delta$ over ground set $[\ell]$ is the ideal
generated by all monomials $\xx_\sigma = \displaystyle{\prod_{i \in
    \sigma}} x_i$ such that $\sigma \not\in \Delta$. 
Conversely, each squarefree monomial ideal in $T$ is the Stanley-Reisner 
ideal of
some simplicial complex on ground set $[\ell]:=\{ 1,
\ldots, \ell\}$.

If $f \in T$ and $\preceq$ is a term
order, then we denote by $\iin_\preceq (f)$ the leading monomial of
$f$ with respect to $\preceq$ (i.e., $0$ if $f = 0$ and the largest
monomial with respect to $\preceq$ occurring in $f$ otherwise). For an
ideal $J$, we write $\iin_\preceq(J)$ for the {\em initial ideal} of
$J$, i.e., the ideal generated by $\{\iin_\preceq(f)
: f \in J\}$. Clearly, 
$\iin_\preceq(J)$ is a monomial ideal. 

If $\Delta$ and $\Gamma$ are simplicial complexes on disjoint ground
sets, then the {\em join} $\Delta * \Gamma$ is the simplicial complex 
$$
\Delta * \Gamma := 
\{ \sigma \cup \tau : \sigma \in \Delta, \tau \in \Gamma \}.
$$ 

We are now in position to state the main results of this paper. 
 
\begin{theorem} \label{maintheorem}
  Let $1 \leq 2r \leq n$. Then there is a term order $\preceq$ on
  the monomials in $S_n$ such that
  $$
  \iin_{\preceq} (\Pfaff_{n,r}) 
  = I_{\Delta_{n,r-1} * 2^{\Theta_{n,r-1}}},
  $$
  where $\Delta_{n,r-1}$ is a simplicial sphere of dimension
  $r(n-2r-1)-1$ and $2^{\Theta_{n,r-1}}$ is the full simplex on a set
  $\Theta_{n,r-1}$ of size $nr$.
\end{theorem}
To describe the simplicial complexes $\Delta_{n,r}$ appearing in the
formulation of Theorem \ref{maintheorem}, we need to recall
some more facts from the theory of simplicial complexes. 

An element $\sigma \in \Delta$ of a simplicial complex is called
a {\em face} of $\Delta$. The {\em dimension} $\dim(\sigma)$ of
$\sigma$ is defined as $\dim(\sigma) := \# \sigma -1$. The dimension
of $\Delta$  is defined as $\dim(\Delta) = \max \{ \dim(\sigma) :
\sigma \in \Delta\}$. If all inclusionwise maximal faces $\sigma$ of
$\Delta$ are of dimension $\dim(\Delta)$, then $\Delta$ is called {\em
  pure}.

Let $\Omega_n = \{ (i,j) : 1 \leq i <j \leq n\}$. We can visualize the
elements of $\Omega_n$ as edges and diagonals in a convex $n$-gon with
vertices labelled in clockwise order $1 , \ldots, n$. 
Clearly, $\Omega_n$ is in bijection with the variables of the polynomial
ring $S_n$. Therefore, for any simplicial complex $\Delta$ on ground
set $\Omega_n$, we may consider $I_\Delta$ as an ideal in $S_n$. 
We introduce a distance function $d$ on $[n]$ by setting
$$
d_{ij} = \min \{|j-i|, |n+i-j|\}.
$$
Equivalently, if we again regard $i$ and $j$ as vertices of a convex
$n$-gon, then $d_{ij}-1$ is the minimum of the number of vertices 
on the left-hand side and the number of vertices on the right-hand
side of the diagonal through $i$ and $j$.
We denote by $\Omega_{n,r}$ the set of elements $(i,j) 
\in \Omega_n$ such that $d_{ij} > r$.

Let $\Delta_{n,r}$ be the simplicial complex on ground set $\Omega_{n,r}$
whose simplices are those subsets $\sigma \subseteq \Omega_{n,r}$ for which
there is no $\tau \subseteq \sigma$ such that $\# \tau = r+1$ and 
any two diagonals in $\tau$ intersect; by intersection we mean
transversal intersection in the interior of the $n$-gon. 

\begin{proposition} \label{sphere}
  Let $0 \leq 2r \leq n-1$.
  \begin{itemize}
    \item[(i)] \cite{DGJM} The simplicial complex $\Delta_{n,r}$
               triangulates a PL-sphere of dimension
               $\dim(\Delta_{n,r}) = r(n-2r-1)-1$.
    \item[(ii)] \cite{J} The number of faces of $\Delta_{n,r}$ of 
                dimension $r(n-2r-1)-1$ is given by  
                $$\prod_{1 \leq i \leq j \leq n-2r-1} \frac{2r+i+j}{i+j}.$$
    \end{itemize}
\end{proposition}

We close this section by recalling some facts that link invariants of
simplicial complexes and invariants of rings.
If $J$ is a finitely generated ideal in a ring $T$ such that $J$ is
homogeneous
with respect to the standard grading on monomials, then $R := T/J \cong 
\displaystyle{\bigoplus_{i \geq 0} R_i}$ is a standard graded ring,
where $R_i$ is the $k$-vectorspace of elements of degree $i$.  
In this situation, the Hilbert-series $\Hilb(R,t) = 
\displaystyle{\sum_{i \geq 0}} \dim_k R_i t^i$ 
is a rational function 
$$\Hilb(R,t) = \frac{h_R(t)}{(1-t)^d},$$ 
where $h_R(t)$ is a polynomial and $d$ is the Krull dimension 
$\dim(R)$ of $R$. The {\em
$a$-invariant} $a(R)$ of $R$ is then
defined as the difference of the degree of the polynomial $h_R(t)$ and
$d$. If $\Delta$ is a simplicial complex, then $\dim(T/I_\Delta) =
\dim(\Delta)+1$ \cite[Theorem 5.1.4]{BH}. The number $h_R(1)$ is called
the multiplicity $e(R)$ of $R$. If $J = I_\Delta$, then
$e(T/I_\Delta)$ is the number of faces $\sigma \in \Delta$ such that
$\dim(\sigma) = \dim(\Delta)$.

\section{Consequences of Herzog and Trung's Term Order}

Let $\Sigma_{n,r}$ be the simplicial complex on ground set
$\Omega_n$ such that the minimal nonfaces of
$\Sigma_{n,r}$ are those sets $\tau$ of size $r+1$ such that 
any two diagonals $ab$ and $cd$ in $\tau$ are nested, meaning that
$a < c < d < b$ (assuming $a<b$, $c<d$, and $a \le c$). 
\begin{proposition}[\cite{HT}] \label{herzogtrung}
  Let $1 \leq 2r \leq n$. Then there is a term order $\preceq'$ on
  the monomials in $S_n$ such that
  $$
  \iin_{\preceq'} (\Pfaff_{n,r}) 
  = I_{\Sigma_{n,r-1}}.
  $$
\end{proposition}
Using the term order in Proposition~\ref{herzogtrung}, Herzog and
Trung \cite{HT} derived a determinant formula for the multiplicity of
$S_n/\Pfaff_{n,r}$. Using the same term order, Ghorpade and
Krattenthaler \cite{GK} 
obtained a similar determinant formula, which they simplified to
the following product formula via a determinant identity due to
Krattenthaler \cite[Th. 7]{K1}:
\begin{proposition}[\cite{GK}] \label{multiplicity}
  We have that
  $$e(S_n/\Pfaff_{n,r}) =
  \prod_{1 \leq i \leq j \leq n-2r+1} \frac{2(r-1)+i+j}{i+j}.$$
\end{proposition}
Earlier, Harris and Tu \cite{HaTu} discovered a different determinant
formula and also a different product formula for the same multiplicity.
 
By Proposition~\ref{multiplicity}, the multiplicities of the two rings
$S_n/\Pfaff_{n,r}$ 
and $S_n/I_{\Delta_{n,r-1} * 2^{\Theta_{n,r-1}}}$ coincide; apply
Proposition~\ref{sphere} (ii).

De Negri \cite{DN1} derived a determinant formula for the Hilbert
series of $S_n/\Pfaff_{n,r}$ and hence for the $h$-vector of
$\Sigma_{n,r-1}$. Ghorpade and Krattenthaler \cite{GK} recovered this
formula and also provided two alternative determinant formulas for the
same Hilbert series.

By Proposition~\ref{herzogtrung}, Theorem~\ref{maintheorem} implies
the following:
\begin{corollary} \label{samehvector}
  Let $1 \leq 2r \leq n$. Then the two simplicial complexes
  $\Delta_{n,r-1}$ and $\Sigma_{n,r-1}$ have the same $h$-vector. 
\end{corollary}

\section{The New Term Order}
\label{termorder}

Using the distance function $d$ from Section \ref{definitions}, 
we define a linear order  
$\preceq$ on the set of variables such that  
$\varx{i}{j} \prec \varx{k}{l}$ whenever $d_{ij} < d_{kl}$
and extend $\preceq$ to monomials using reverse lexicographic
term order.
Precisely, $m \prec m'$ if and only if either $\deg m < \deg m'$
or $\deg m = \deg m'$ and $m$ is lexicographically smaller
than $m'$,
where the variables in $m$ and $m'$ are arranged in {\em increasing}
order. Note that the term order $\preceq$ does not depend on
the parameter $r$.

For example, for $n=5$, we may choose the order 
\[
\varx{1}{2} \prec \varx{2}{3} \prec \varx{3}{4} \prec \varx{4}{5} \prec
\varx{1}{5} \prec \varx{1}{3} \prec \varx{2}{4} \prec \varx{3}{5} \prec
\varx{1}{4} \prec \varx{2}{5}.
\]
The order of the terms that appear in the Pfaffian of some 
$\Amatrix_{i_1, i_2, i_3, i_4}$ then becomes
\begin{eqnarray*}
& & \varx{1}{2} \varx{3}{4} \prec \varx{1}{2} \varx{4}{5} \prec 
\varx{1}{2} \varx{3}{5} \prec \varx{2}{3} \varx{4}{5} \prec 
\varx{2}{3} \varx{1}{5} \\
&\prec\!\!\!&
\varx{2}{3} \varx{1}{4} \prec
\varx{3}{4} \varx{1}{5} \prec \varx{3}{4} \varx{2}{5} \prec 
\varx{4}{5} \varx{1}{3} \prec \varx{1}{5} \varx{2}{4}  \\
&\prec\!\!\!&
\varx{1}{3} \varx{2}{4} \prec \varx{1}{3} \varx{2}{5} \prec
\varx{2}{4} \varx{3}{5} \prec \varx{3}{5} \varx{1}{4} \prec 
\varx{1}{4} \varx{2}{5}.
\end{eqnarray*}

Note that the monomials in the last row correspond to the minimal
nonfaces of $\Delta_{5,1}$.

\begin{lemma} \label{termorder:lem}
  Let $1 \leq r \leq \frac{n}{2}$ and set $\Theta_{n,r} := 
  \Omega_n \setminus \Omega_{n,r}$.
  Then 
  $I_{\Delta_{n,r-1} * \Theta_{n,r-1}} \subseteq \iin_\preceq(\Pfaff_{n,r})$.
\end{lemma}
\begin{proof}
  The lemma is trivially true for $r=1$; thus assume that $r
  > 1$.  
  For every minimal nonface $\sigma$ of $\Delta_{n,r-1}$, we need to
  prove that $m = \prod_{ij \in \sigma} \varx{i}{j}$ is the leading term of
  some element in $\Pfaff_{n,r}$. Write 
  \[
  \sigma = \{j_1j_{r+1}, j_2 j_{r+2}, \ldots, j_rj_{2r}\},
  \]
  where
  $j_1 < j_2 < \ldots < j_{2r}$.
  It suffices to prove that $m$
  is maximal among all terms in the Pfaffian of $\Amatrix_{j_1, j_2,
  \ldots, j_{2r}}$.
  Assume to the contrary that there is a term $m'$ in the given
  Pfaffian such that $m \prec m'$. 

  By symmetry, we may assume that 
  $\varx{j_1}{j_{r+1}}$ is minimal among the variables in
  $m$ and that $j_{r+1}-j_1 \le  n+j_1- j_{r+1}$.
  The latter property implies that $\varx{j_k}{j_l} \preceq
  \varx{j_1}{j_{r+1}}$ 
  whenever $1 \le k < l \le r+1$.
  By assumption, $m'$ contains no such term $\varx{j_k}{j_l}$,
  except possibly $\varx{j_1}{j_{r+1}}$ itself.
    
  Recall that we may identify any term in the Pfaffian of 
  $\Amatrix_{j_1, j_2, \ldots, j_{2r}}$ with a perfect matching on the set 
  $\{j_1, \ldots, j_{2r}\}$. We say that $j_k$ is {\em matched} with
  $j_l$ in such a term if the term contains the variable
  $\varx{j_k}{j_l}$.
    
  We identify two cases:
  \begin{itemize}
    \item
      $m'$ does not contain $\varx{j_{1}}{j_{r+1}}$. 
      Then each $i \in [r+1]$ must be matched with some element from
      $[r+2,2r]$, because otherwise $m'$ would contain variables that
      are strictly smaller than $\varx{j_{1}}{j_{r+1}}$. However, this is a
      contradiction, because $[r+2,2r]$ has size only $r-1$.
    \item
      $m'$ does contain $\varx{j_{1}}{j_{r+1}}$. Consider the two
      monomials $\hat{m} = m/\varx{j_{1}}{j_{r+1}}$ and $\hat{m'} =
      m'/\varx{j_{1}}{j_{r+1}}$, which both appear in the Pfaffian of 
      $\Amatrix_{j_2, \ldots, j_r, j_{r+2}, \ldots, j_{2r}}$. By
      induction on $r$, we have that 
      $\hat{m'} \preceq \hat{m}$, because 
      $\sigma - j_1j_{r+1}$ is a minimal nonface of $\Delta_{n,r-2}$.
      Yet, this implies that
      \[
      m' = \varx{j_{1}}{j_{r+1}}\hat{m'} \preceq 
      \varx{j_{1}}{j_{r+1}}\hat{m} = m,
      \]
      which is a contradiction.
  \end{itemize}
\end{proof}

\section{Proof of the Main Theorem} \label{proof}

To prove Theorem \ref{maintheorem}, we
need the following lemma about the inclusion of monomial ideals.
A more general argument of this type appears in \cite[Lemma 4.2]{CHT}.

\begin{lemma} \label{monomialinclusion}
  Let $T = k[x_1, \ldots, x_\ell]$ be the polynomial ring in $\ell$
  variables. Suppose that $I \subseteq J$ are monomials ideals in $T$
  such that the following hold:
    \begin{itemize}
      \item[(i)] $\dim(T/I) = \dim(T/J)$.
      \item[(ii)] $e(T/I) = e(T/J)$.
      \item[(iii)] $I = I_\Delta$ for a pure simplicial complex $\Delta$
                   on ground set $[\ell]$. 
    \end{itemize}

    Then $I = J$.
\end{lemma}
\begin{proof}
  Since $\dim(T/I_\Delta) = \dim(\Delta) +1$, it follows by (i) that
  $\dim(\Delta) = \dim(T/J) -1$. Let $J^\pol$ be the polarization of 
  $J$ that lives in the polynomial ring $T'$ in $\ell'$ variables. We refer
  the reader to \cite{F1} for the definition and basics of polarization.
  We can assume that $T \subseteq T'$. 
  It is well known \cite{F1} that 
  $\dim(T'/J^\pol) = \dim(T/J) + (\ell' - \ell)$.
  Let $I'$ be the ideal generated by $I$ in $T'$. Then $I'$ is the 
  Stanley-Reisner ideal of a simplicial complex 
  $\Delta' = \Delta * 2^\Theta$ for some set $\Theta$ of cardinality
  $\ell' - \ell$. In particular, 
  $$\dim(\Delta') = \dim(\Delta) + (\ell' - \ell)
  = \dim(T'/J^\pol) - 1.$$
  Since $J^\pol$ is a squarefree monomial  
  ideal, there is a simplicial complex $\Gamma$ on ground set
  $[\ell']$ for which $J^\pol = I_\Gamma$.
  By $I' = I_{\Delta'} \subseteq J^\pol = I_\Gamma$, it follows
  that $\Gamma \subseteq \Delta'$. Moreover, 
  $\dim(\Gamma) = \dim(\Delta')$.

  The number of faces of $\Delta$ in top dimension 
  $\dim(\Delta)$
  coincides
  with the number of faces of 
  $\Delta'$ in top dimension 
  $\dim(\Delta')$
  and is equal to the
  multiplicity $e(T/I) = e(T'/I')$ of $T/I$ and $T'/I'$. Since the
  multiplicity does not change under polarization \cite{F1}, it
  follows that $e(T/J) = e(T'/J^\pol)$. 
  Since $e(T'/J^\pol)$ equals the number of faces of $\Gamma$ in
  top dimension $\dim(\Gamma) = \dim(\Delta')$, we may conclude by
  (i) and (ii) that $\Gamma$ and $\Delta'$ have the same number of
  faces in this dimension. By $\Gamma \subseteq \Delta'$, it follows
  that these faces are indeed the same. From the fact (iii) that
  $\Delta$ is pure, 
  we obtain that any face of $\Delta$ is contained in a
  top-dimensional face. Since purity of $\Delta$ implies 
  purity of $\Delta'$, the same holds for $\Delta'$.
  But then $\Gamma \subseteq \Delta'$ implies $\Gamma =
  \Delta'$. Therefore, the minimal monomial generating set of $J^\pol
  = I_\Gamma$ uses only the variables from $T$. Thus $J = J^\pol \cap
  T = I_\Delta = I$.
\end{proof}

\begin{proof}[Proof of Theorem \ref{maintheorem}:]
  Let $S_n$, $\Pfaff_{n,r+1}$, $\Delta_{n,r}$, $\Omega_n$, $\Omega_{n,r}$,
  $\Theta_{n,r}$ and $\preceq$ be as in Sections \ref{definitions}
  and \ref{termorder}.

  We set $I:=I_{\Delta_{n,r} * \Theta_{n,r}}$ and $J :=
  \iin_{\preceq}(\Pfaff_{n,r+1})$. 
 
  \begin{itemize}
    \item[$\triangleright$] By Lemma \ref{termorder:lem}, we know that
    $I \subseteq J$. 
    \item[$\triangleright$] By Proposition \ref{sphere} (i), we know
         that $\Delta_{n,r}$ and hence $\Delta_{n,r} *
         \Theta_{n,r}$ are pure simplicial complexes. Moreover, the
         dimension of $\Delta_{n,r}$ is $r(n-2r-1)-1$, which yields
         that 
	 $$\dim(\Delta_{n,r} * 2^{\Theta_{n,r}})  
         = r(n-2r-1) + \#\Theta_{n,r} - 1.$$ 
	 Simple enumeration shows that $\# \Theta_{n,r} = nr$; 
         thus 
	 $$\dim(\Delta_{n,r} * 2^{\Theta_{n,r}}) =r(2n-2r-1) - 1.$$ 
	 Hence
         $\dim(S_n/I) = r(2n-2r-1)$, which coincides \cite{HT} with
	 $\dim(S_n/\Pfaff_{n,r+1})$.  
         Since by general Gr\"obner basis theory
	 $$\dim(S_n/\Pfaff_{n,r+1}) = 
         \dim(S_n/\iin_\preceq (\Pfaff_{n,r+1})) = \dim(S_n/J),$$ 
	 it follows that $\dim(S_n/I) = \dim(S_n/J)$. 
    \item[$\triangleright$] The number of maximal faces of $\Delta_{n,r}$ 
         coincides with the number of maximal faces of $\Delta_{n,r} * 2^{\Theta_{n,r}}$
         and therefore with $e(S_n/I)$. Again general Gr\"obner basis theory says that
         $$e(S_n/\Pfaff_{n,r+1}) = e(S_n/\iin_\preceq
         (\Pfaff_{n,r+1})) = e(S_n/J).$$
         By Proposition \ref{sphere} (ii) and Proposition
         \ref{multiplicity},
	 it follows that $e(S_n/I) = e(S_n/J)$. 
  \end{itemize}
  Therefore, the ideals $I$ and $J$ satisfy the assumptions of Lemma 
  \ref{monomialinclusion}. Thus $I = J$, which concludes the proof.
\end{proof}
Recently, in a more general setting, Krattenthaler \cite{K2} obtained a
degree-preserving bijection between the sets of monomials in the two
rings $S_n/I_{\Delta_{n,r-1} * 2^{\Theta_{n,r-1}}}$ and
$S_n/I_{\Sigma_{n,r-1}}$. The bijection
is based on a variation of the Robinson-Schensted-Knuth
correspondence. Using this bijection, one may prove
Theorem~\ref{maintheorem} without using Propositions~\ref{sphere} (ii) 
and \ref{multiplicity}. 

\section{Determinantal Ideals}
\label{determinantal}

Using Lemma~\ref{monomialinclusion}, we prove that 
certain simplicial complexes related to a given determinantal ideal
have the same $h$-vector, thereby generalizing enumerative results due
to the first author \cite{J} presented in Proposition~\ref{content}
below. This time, the derived initial ideals are not spherical in
general.

Let
$$M = 
\left( 
\begin{array}{ccccc} 
\varx{1}{1}    & \varx{1}{2}   & \varx{1}{3}    & \cdots \\
\varx{2}{1}   & \varx{2}{2}    & \varx{2}{3}    & \cdots \\               
\varx{3}{1}   & \varx{3}{2}  &   \varx{3}{3}    & \cdots \\
\cdots    & \cdots   & \cdots    & \cdots \\
\cdots    & \cdots   & \cdots    & \cdots \\
\end{array}
\right)$$
be a generic matrix indexed by $\mathbb{P}^2$, where $\mathbb{P} =
\{1,2, 3, \ldots\}$. 
Let $\Lambda$ be a finite subset of $\mathbb{P}^2$
such that the following hold:
\begin{itemize}
\item
  If $ab := (a,b) \in \Lambda$, then $cb \in \Lambda$ for $1 \le c \le
  a$.
\item
  If $ab_1,ab_2 \in \Lambda$ and $b_1\leq b_2$, then $ad \in
  \Lambda$ for $b_1 \leq d \leq b_2$.
\end{itemize}
We refer to $\Lambda$ as a {\em stack polyomino}. 
If in addition $a1 \in \Lambda$ whenever $ab \in \Lambda$ for some
$b$, then $\Lambda$ is a {\em Ferrers diagram}.

Define $M(\Lambda) = (m_{ab})$ by
$$m_{ab} =
\left\{ 
\begin{array}{cl}
  \varx{a}{b} & \mbox{if $ab \in \Lambda$;}\\
  0 & \mbox{otherwise.}
\end{array}
\right.
$$
For $r \ge 1$ and index sets
$\alpha = \{a_1, \ldots, a_{r}\}$ and $\beta = \{b_1, \ldots, b_r\}$
of size
$r$, we denote by 
$M_{\alpha,\beta}(\Lambda)$ the submatrix of
$M(\Lambda)$ constructed by selecting the rows indexed by $\alpha$ 
and the columns indexed by $\beta$.  
Define
$$
D_{\alpha,\beta}(\Lambda) =
\left\{ 
\begin{array}{cl}
  \det M_{\alpha,\beta}(\Lambda) & \mbox{if
  $a_i b_j \in \Lambda$ for $1\le i,j \le r$;}\\ 
  0 & \mbox{otherwise.}
\end{array}
\right.
$$
Let $R_{\Lambda} = k[\varx{a}{b} : ab \in \Lambda]$ be
the polynomial ring over the field $k$ in the variables
indexed by $\Lambda$.
We denote by $\Determ_{\Lambda,r}$ the ideal in $R_\Lambda$ generated by
the polynomials $D_{\alpha,\beta}(\Lambda)$ for all possible choices of
index sets $\alpha$ and $\beta$ of size $r$. Note that all but
finitely many of these polynomials are zero.

Two elements $ab$ and $cd$ form a 2-{\em diagonal} in 
$\Lambda$ if $a<c$, $b<d$, and the $2 \times 2$ square
$\{ab, ad, cb, cd\} = \{a,c\} \times \{b,d\}$ is a subset
of $\Lambda$. More generally,
$a_1 b_1, \ldots, a_rb_r$ form an {\em $r$-diagonal}
in $\Lambda$ if the following hold:
\begin{itemize}
\item
  $a_1 < a_2 < \ldots < a_r$ and $b_1 < b_2 < \ldots < b_r$
\item
  The $r \times r$ square
  $\{ a_ib_j : i, j \in [1,r] \} =
  \{ a_1, \ldots, a_r \} \times
  \{ b_1, \ldots, b_r \}$
  is a subset of $\Lambda$.
\end{itemize}
For $r \ge 1$, define $\Delta_{\Lambda,r-1}$ as the family of
subsets $\sigma$ of $\Lambda$ such that $\sigma$ does not contain
any set forming an $r$-diagonal in $\Lambda$. Clearly,
$\Delta_{\Lambda,r-1}$ is a simplicial complex.

Define a linear order on the set of variables in $\Lambda$ such that
$\varx{i}{j} \prec \varx{k}{l}$ whenever $i>k$ and also whenever $i=k$
and $j<l$. Extend $\preceq$ to monomials using reverse
lexicographic term order in the same manner as in
Section~\ref{termorder}. Again, the term order $\preceq$ does not
depend on $r$.
\begin{lemma} \label{termorder2:lem}
  Let $r \geq 1$.
  Then 
  $I_{\Delta_{\Lambda,r-1}} \subseteq
  \iin_\preceq(\Determ_{\Lambda,r})$.
\end{lemma}
\begin{proof}
  The lemma is trivially true for $r=1$; thus assume that $r
  > 1$.  
  As in the proof of Lemma~\ref{termorder:lem}, 
  we need to prove that $m = \prod_{ab \in \sigma} \varx{a}{b}$ is the
  leading term of some element in $\Determ_{n,r}$
  for every minimal nonface $\sigma$ of $\Delta_{\Lambda,r-1}$.
  Write 
  \[
  \sigma = \{a_1b_{1}, a_2 b_2, \ldots, a_rb_r\},
  \]
  where
  $a_1 < a_2 < \ldots < a_r$ and $b_1 < b_2 < \ldots < b_r$.
  Let $\alpha = \{a_1, \ldots, a_r\}$ and 
  $\beta = \{b_1, \ldots, b_r\}$.
  It suffices to prove that $m$ is maximal among all terms in
  $D_{\alpha,\beta}(\Lambda)$.
  Now, each monomial in $D_{\alpha,\beta}(\Lambda)$ contains one
  variable from each row in $\alpha$ and 
  one variable from each column in $\beta$. As a consequence, 
  the smallest variable in any term of $D_{\alpha,\beta}(\Lambda)$
  appears in row $a_r$. In particular, the maximal term of
  $D_{\alpha,\beta}(\Lambda)$ contains the element
  $\varx{a_r}{b_r}$. Proceeding by induction on $r$ with $D_{\alpha
  \setminus \{a_r\},\beta \setminus \{b_r\}}(\Lambda)$ in
  the same manner as in the proof of Lemma~\ref{termorder2:lem}, we
  conclude that $m$ is indeed maximal in $D_{\alpha,\beta}(\Lambda)$.
\end{proof}

We may identify any given stack polyomino $\Lambda$ with the sequence
 $(\lambda_1,\lambda_2, \lambda_3, \ldots)$ defined by $\lambda_j =
 \max \{ i : (i,j) \in \Lambda\}$. Note that $\lambda_j$ is the number
 of elements in column $j$ in $\Lambda$. 
 The {\em content} of $\Lambda$ is the Ferrers diagram obtained by
 arranging the elements $\lambda_1, \lambda_2, \lambda_3, \ldots$ in
 decreasing order. 

\begin{proposition}[Jonsson \cite{J}] \label{content}
  Let $\Lambda$ be a stack polyomino and let $\Lambda'$ be its
  content. Then $\Delta_{\Lambda,r-1}$ and $\Delta_{\Lambda',r-1}$ are
  pure complexes of the same dimension with the same number of maximal
  faces. 
\end{proposition}
\begin{proposition}[Herzog and Trung \cite{HT}] \label{ferrers}
  Let $\Lambda$ be a Ferrers diagram. Then
  $I_{\Delta_{\Lambda,r-1}} = \iin_\preceq(\Determ_{\Lambda,r})$.
\end{proposition}
\begin{theorem} \label{determ:thm}
  Let $\Lambda$ be a stack polyomino. Then
  $I_{\Delta_{\Lambda,r-1}} = \iin_\preceq(\Determ_{\Lambda,r})$.
\end{theorem}
\begin{proof}
  Set $I:= I_{\Delta_{\Lambda,r-1}}$ and $J :=
  \iin_\preceq(\Determ_{\Lambda,r})$.
  \begin{itemize}
    \item[$\triangleright$] By Lemma \ref{termorder2:lem}, we know
    that $I \subseteq J$. 
    \item[$\triangleright$] By Proposition \ref{content}, we have that
      $\dim(R_\Lambda/I) = \dim(R_{\Lambda'}/{I'})$, where
      $\Lambda'$ is the content of $\Lambda$ and $I' :=
      I_{\Delta_{\Lambda',r-1}}$.
      By Proposition~\ref{ferrers}, it follows that
      $$\dim(R_{\Lambda'}/{I'}) =
      \dim(R_{\Lambda'}/\Determ_{\Lambda',r+1}).$$  
      Since $R_\Lambda/\Determ_{\Lambda,r+1}$ and 
      $R_{\Lambda'}/\Determ_{\Lambda',r+1}$ are isomorphic,
      we conclude that this equals 
      $$\dim(R_{\Lambda}/\Determ_{\Lambda,r+1}) = \dim(R_\Lambda/J);
      $$
      hence $\dim(R_\Lambda/I) = \dim(R_\Lambda/J)$.
    \item[$\triangleright$] By 
      Proposition \ref{content}, we have that
      $e(R_\Lambda/I) = e(R_{\Lambda'}/{I'})$.
      Another application of Proposition~\ref{ferrers} yields that
      $$e(R_{\Lambda'}/{I'}) = 
      e(R_{\Lambda'}/\Determ_{\Lambda',r+1}) =
      e(R_{\Lambda}/\Determ_{\Lambda,r+1}) = e(R_\Lambda/J). 
      $$
      Summarizing, we deduce that $e(R_\Lambda/I) = e(R_\Lambda/J)$.
  \end{itemize}
  Applying Lemma \ref{monomialinclusion}, we obtain that $I = J$, and
  we are done.
\end{proof}
\begin{corollary} \label{determ:cor}
  Let $\Lambda$ and $\Lambda'$ be stack polyominoes with the same
  content. Then $\Delta_{\Lambda,r-1}$ and $\Delta_{\Lambda',r-1}$
  have the same $f$-vector.
\end{corollary}
The bijection of Krattenthaler \cite{K2} mentioned at the end of
Section~\ref{proof} applies also to this situation for the special
case that $\Lambda$ is a Ferrers diagram and 
$\Lambda'$ is the (horizontal) reflection of $\Lambda$. 
In particular, Corollary~\ref{determ:cor} is a consequence of this
bijection for the given special case.

\end{document}